\documentclass[11pt,leqno]{article} 

\newtheorem{thm}{Theorem}[section]
\newtheorem{lma}{Lemma}[section]

\newcommand{\beqa}{\begin{eqnarray}}
\newcommand{\eeqa}{\end{eqnarray}}

\newcommand{\pf}{\noindent {\bf Proof:} $\s$ }
\newcommand{\epf}{ \hfill$\diamondsuit$ \medskip}

\newcommand{\beq}{\begin{equation}}
\newcommand{\eeq}{\end{equation}}
\newcommand{\lbl}{\label}
\newcommand{\s}{\; \;}

\newcommand{\la}{\lambda}

\newcommand{\ra}{\rightarrow}

\newcommand{\p}{\varphi}

\usepackage{graphics}

\title{A global solution curve for a class of periodic problems, including the relativistic pendulum }

\author{
Philip Korman   \\ 
Department of Mathematical Sciences \\ 
University of Cincinnati \\ 
Cincinnati Ohio 45221-0025 \\
}

\date{}

\begin{document}

\maketitle
\begin{abstract} 
Using continuation methods, we study  the global solution structure of periodic solutions for a class of periodically forced  equations, generalizing the case of relativistic pendulum. We obtain results on the existence and multiplicity of periodic solutions. Our approach is suitable for numerical computations, and in fact we present some numerically computed bifurcation diagrams illustrating our results. 
 \end{abstract}

\begin{flushleft}
\end{flushleft}

\begin{flushleft}
\end{flushleft}

\section{Introduction}
\setcounter{equation}{0}
\setcounter{thm}{0}
\setcounter{lma}{0}

There is a considerable  recent interest in periodic solutions of a class of equations generalizing the relativistic pendulum  equation
\beqa
\lbl{i1}
& (\p (u'))'+\la f(u)u'+ kg(u)=h(t)=\mu + e(t), \\
& u(t+T)=u(t), \s  u'(t+T)=u'(t)\,. \nonumber
\eeqa
For relativistic pendulum we have $\p (z)=\frac{z}{\sqrt{1-z^2}}$, $g(u)=\sin u$, and $f(u)=1$.
More generally, the function $\p (z) : (-a,a) \ra R$ is assumed to be increasing, of class $C^1$, and $\lim _{z \ra \pm a} \p(z)=\pm \infty$. Without restricting the generality, we shall also assume that $\p (0)=0$.  The function $h(t) \in C(R)$ is assumed to be $T$-periodic, $f(u) \in C(R)$, the friction $\la$ is constant. We decompose $h(t)=\mu + e(t)$, with $\mu = \frac1T \int_0^T h(t) \, dt$, and $\int_0^T e(t) \, dt=0$.
P.J. Torres \cite{T1} proved existence of at least two $T$-periodic solutions, assuming that $aT<2 \sqrt3 $ and $ |\mu|<k \left( 1-\frac{aT}{2 \sqrt3} \right)$. Then
C. Bereanu, P.  Jebelean and J. Mawhin \cite{B} have improved these conditions to read: $aT<\pi \sqrt3$, $ |\mu|<k \cos \frac{aT}{2  \sqrt3}$.
H. Brezis and J. Mawhin \cite{B1} have proved existence of at least one solution,  under different  conditions, and their result was recently extended by C. Bereanu  and P.J.  Torres \cite{BT}.
Our first result is  a more elementary proof of the above result of C. Bereanu, P.  Jebelean and J. Mawhin \cite{B}, which uses only Schauder's fixed point theorem.
\medskip

We seek to understand the shape of the global solution curve, to explain why some restrictions on $\mu$ are necessary for the existence of solutions, and how multiple solutions are connected. Let us decompose the solution $u(t)=\xi + U(t)$, with $\xi = \frac1T \int_0^T u(t) \, dt$, and $\int_0^T U(t) \, dt=0$. We study the global solution curve for the problem (\ref{i1}), i.e., $\mu=\mu(\xi)$. We give conditions under which $\xi$ is a {\em global parameter}, which means that for each  $\xi \in (-\infty,\infty)$ there is a {\em unique} pair $(\mu,u(t))$ solving (\ref{i1}). To establish  that, we continue solutions of (\ref{i1}) back in $k$ on curves of fixed average $\xi$, similarly to author's recent papers \cite{K1}, \cite{K} and \cite{K2}, and at $k=0$ we have a complete description of $T$-periodic solutions,  in particular we have the existence and uniqueness of $T$ periodic solutions of any average, which implies the uniqueness of $(\mu,u(t))$.
We then study properties of the curve $\mu=\mu(\xi)$, depending on $g(u)$. If $g(u)$ is a periodic function, it turns out that so is $\mu(\xi)$, with the same period. If  $g(u)$ tends to finite limits, as $u \ra \pm \infty$, then $\mu(\xi)$ tends to finite limits, as $\xi \ra \pm \infty$. We illustrate our results by numerical computations, and we discuss in detail their implementation. Our computations show that solutions of (\ref{i1}) have very small variation, i.e., they are close to constant solutions of the algebraic equation $kg(u)=\mu$.

\section{Preliminary results}
\setcounter{equation}{0}
\setcounter{thm}{0}
\setcounter{lma}{0}

We record the following simple observation.
\begin{lma}\lbl{lma:1}
Let $e(t)$ be a  given continuous function of period $T$. Then the function $\int _0^t e(s) \, ds$ is $T$-periodic if and only if $\int_0^T e(s) \, ds=0$.
\end{lma}

The next lemma deals with periodic solutions of ``$\p$-linear" equations.
\begin{lma}\lbl{lma:2}
Consider the  problem 
\beq
\lbl{1}
(\p (u'(t)))'+\la u'(t)=e(t),
\eeq
where $\p : (-a,a) \ra R$ is an increasing homeomorphism, of class $C^1$, with $\p (0)=0$, and 
 $e(t)$ is a given continuous function of period $T$, of zero average, i.e.,  $\int_0^T e(s) \, ds=0$, and $\la \geq 0$, a constant. Then the problem (\ref{1}) has a family of  $T$-periodic solutions of the form $u=u_0(t)+c$, where $c$ is any constant. This family exhausts the set of $T$-periodic solutions. In particular, one can select a  unique $T$-periodic solution of any average.
\end{lma}

\pf
{\bf Case 1}. $\la =0$. From (\ref{1})
\[
u'(t)=\p ^{-1} \left[\int_0^t e(s) \, ds+C \right], \s \mbox{$C$ is an arbitrary constant} \,.
\]
Observe that $u'(t)$ is a $T$-periodic function, whose values belong to the domain of $\p$. We can choose $C=C_0$, such that 
\[
\int_0^T \p ^{-1} \left[\int_0^t e(s) \, ds+C_0 \right] \, dt=0 \,.
\]
(This integral is positive (negative) for $|C|$ large and $C$ positive (negative). Existence of $C_0$ follows by continuity.) By monotonicity of $\p ^{-1}$, $C_0$ is unique. Then
\[
u(t)=\int_0^t \p ^{-1} \left[\int_0^t e(s) \, ds+C_0 \right] \, dt+c
\]
gives us the desired $T$-periodic solutions.
\medskip

\noindent
{\bf Case 2}. $\la >0$. Letting $u'=z$, we rewrite (\ref{1}) as
\beq
\lbl{2} 
\p (z)'+\la z=e(t) \,.
\eeq
It suffices to show that (\ref{2}) has a $T$-periodic solution. Indeed, integrating (\ref{2}), we see that this solution satisfies
\[
\int _0^T z(t) \, dt=0 \,,
\]
and then $u(t)=\int_0^t z(s) \, ds+c$ gives us the desired $T$-periodic solutions. Let $\p (z(t))=p(t)$, i.e., $z=\p ^{-1} (p) \equiv \psi (p)$. Observe that again the values of $z=u'(t)$ belong to the domain of $\p$, and that $\psi (p)$ is a bounded function, which is positive (negative) for $p$ positive (negative).
We rewrite (\ref{2}) as 
\beq
\lbl{3}
p'+\la \psi (p)=e(t) \,.
\eeq
Let us solve (\ref{3}), with $p(0)=p_0$, and call the solution $p(t,p_0)$. If $p_0>0$ and large, then $p(t)>0$ for all $t \in [0,T]$. Then, integrating (\ref{3}),
\[
p(T)<p(0)=p_0 \,.
\]
We conclude that the Poincare map $p_0 \ra p(T,p_0)$ takes the interval $(-p_0,p_0)$ into itself, for $p_0>0$  large. Hence (\ref{3}), and therefore (\ref{2}), have periodic solutions. Since $\psi (p)$ is monotone, (\ref{3}) has a unique $T$-periodic solution.
\epf

The following lemma is known. We present its proof for completeness.

\begin{lma}\lbl{lma:3}
Let $u(t) \in W^{1,2}(R)$ be a $T$-periodic function, with $\int_0^T u(t) \, dt=0$. Then
\beq
\lbl{4}
||u||^2_{L^{\infty}(R)} \leq \frac{T}{12} \int_0^T {u'}^2(t) \, dt \,.
\eeq
 \end{lma}

\pf
Represent $u(t)$ by its Fourier series on $(-T/2,T/2)$ (with $a_0=0$)
\[
u(t)=\Sigma _{n=1}^{\infty} a_n \cos \frac{2n \pi}{T} t+b_n \sin \frac{2n \pi}{T} t \,,
\]
and then 
\beq
\lbl{5}
\int_0^T {u'}^2(t) \, dt=\frac{2 \pi ^2}{T} \, \Sigma _{n=1}^{\infty} (a^2_n+b^2_n)n^2  \,.
\eeq
Applying the Schwarz inequality to the scalar product of the vectors $(a_n,b_n)$ and $(\cos \frac{2n \pi}{T} t,\sin \frac{2n \pi}{T} t)$, we have
\[
||u||_{L^{\infty}(R)} \leq \Sigma _{n=1}^{\infty} \sqrt{a^2_n+b^2_n}=\Sigma _{n=1}^{\infty} \frac{1}{n} n \sqrt{a^2_n+b^2_n} \,.
\]
Using (\ref{5}), we then have
\[
||u||^2_{L^{\infty}(R)} \leq \Sigma _{n=1}^{\infty} \frac{1}{n^2} \, \Sigma _{n=1}^{\infty} (a^2_n+b^2_n)n^2=\frac{\pi ^2}{6} \frac{T}{2 \pi ^2}\int_0^T {u'}^2(t) \, dt \,,
\]
and the proof follows.
\epf

\noindent
{\bf Remark} An even shorter proof can be given by using complex Fourier series. Representing $u(t)=\Sigma _{j \ne 0} c_j e^{\frac{2 \pi}{T} i jt}$ (since $c_0=0$), we have
\[
\| u \| _{L^{\infty}} \leq \Sigma _{j \ne 0} |c_j| \leq \left(\Sigma _{j \ne 0} \frac{1}{j^2} \right)^{1/2} \left(\Sigma _{j \ne 0} j^2 |c_j|^2\right)^{1/2} =  \frac{\pi}{\sqrt 3} \frac{\sqrt T}{2 \pi} \| u'\|  _{L^2} \,. 
\]

We consider classical solutions of the problem
\beq
\lbl{6}
(\p (u'))'+f(u)u'+ \sin u=h(t) \,.
\eeq
The function $\p : (-a,a) \ra R$ is assumed to be increasing of class $C^1$, and $\lim _{u \ra \pm a} \p(u)=\pm \infty$. Without restricting the generality, we shall also assume that $\p (0)=0$.
\medskip

We observe the following simple lemma.
\begin{lma}\lbl{lma:2.3}
Let $u(t)$ be a $T$-periodic solution of 
\[
(\p (u'))'=g(t) \,,
\]
where $g(t) \in C(R)$ is a given $T$-periodic function of zero average. Then there is  a constant $\alpha$, $0<\alpha<a$, such that 
\[
|u'(t)|< \alpha, \s \mbox{for all $t$} \,.
\]
\end{lma}

\pf
Integrating the equation, between any critical point $t_0$ of $u(t)$, and any point $t$, with  $t_0, t \in (0,T)$,
\[
\p (u'(t))=\int _{t_0}^t g(t) \, dt \leq \int _{0}^T |g(t)| \, dt \,.
\]
Hence $u'(t)$ cannot get near $\pm a$ on $(0,T)$, and by periodicity, for all $t$.
\epf

This lemma shows that when one continues the solutions of (\ref{6}), $u'(t)$ stays away from the values where $\p (u')$ is not defined.

The following lemma is known as Wirtinger's inequality. Its proof follows easily by using the complex Fourier series, and the orthogonality of the functions $\{e^{i \omega n t} \}$ on the interval $(0,T)$.
\begin{lma}\lbl{lma:4}
Assume that $f(t)$ is a continuously differentiable function  of period $T$, and of zero average, i.e.  $\int_0^T f(s) \, ds=0$. Then, denoting $\omega =\frac{2\pi}{T}$, 
\[
\int _0^T {f'}^2(t) \, dt \geq \omega ^2 \int _0^T f^2(t) \, dt.
\]
\end{lma}

We shall denote by $W_T^{1,2}$ the subset of $W^{1,2}(R)$, consisting of $T$-periodic functions, with the norm $||u||^2_{W_T^{1,2}}=\int_0^T (u^2+{u'}^2) \, dt$.

\begin{lma}\lbl{lma:3.3}
Let $u(t)$ be a $T$-periodic solution of zero average of 
\[
(\p (u'))'=g(t) \,,
\]
where $g(t) \in C(R)$ is a given $T$-periodic function of zero average. Then there is a constant $c_0$ independent of $g(t)$, so that 
\[
||u||_{W_T^{1,2}} \leq c_0 \,.
\]
\end{lma}

\pf
Since $|u'(t)| < a$ for any solution, we have a bound on $\int _0^T {u'}^2 \, dt$, and by Wirtinger's inequality we have a bound on $\int _0^T {u}^2 \, dt$.
\epf

Observe that this is better than what one has for a linear equation $u''=g(t)$, where the bound on $||u||_{W_T^{1,2}}$ does depend on $g(t)$.

\section{Existence of at least two solutions}
\setcounter{equation}{0}
\setcounter{thm}{0}
\setcounter{lma}{0}

We consider classical solutions of the periodic problem 
\beqa
\lbl{7}
& (\p (u'))'+f(u)u'+ k\sin u=h(t)=\mu + e(t), \\
& u(t+T)=u(t), \s  u'(t+T)=u'(t)\,. \nonumber
\eeqa
The function $\p : (-a,a) \ra R$ is assumed to be increasing of class $C^1$, and $\lim _{u \ra \pm a} \p(u)=\pm \infty$. Without restricting the generality, we shall also assume that $\p (0)=0$.  The function $h(t) \in C(R)$ is assumed to be $T$-periodic, $f(u) \in C(R)$. We decompose $h(t)=\mu + e(t)$, with $\mu = \frac1T \int_0^T h(t) \, dt$, and $\int_0^T e(t) \, dt=0$. We present next a simple proof 
of the following result of C. Bereanu, P.  Jebelean and J. Mawhin \cite{B}, see also P.J. Torres \cite{T1}.

\begin{thm} (\cite{B})
Assume that 
\[
aT<\pi \sqrt3, \s |\mu|<k \cos \frac{aT}{2  \sqrt3}.
\]
Then the problem (\ref{7}) has at least two $T$-periodic solutions.
\end{thm}

\pf
Decompose $u(t)=\xi+U(t)$, with $\xi = \frac1T \int_0^T u(t) \, dt$, and $\int_0^T U(t) \, dt=0$. Integrating (\ref{7})
\beq
\lbl{8}
\mu =\frac{k}{T} \int_0^T \sin \left( \xi+U(t) \right) \, dt \,.
\eeq
Using this in (\ref{7})
\beq
\lbl{9}
\s\s\s\s (\p (u'))'+f(u)u'=-k\sin \left( \xi+U(t) \right)+\frac{k}{T} \int_0^T \sin \left( \xi+U(t) \right) \, dt+e(t) \,.
\eeq
The system of (\ref{8}) and (\ref{9}) is equivalent to (\ref{7}), in fact it gives the classical Lyapunov-Schmidt decomposition of (\ref{7}). Let $\overline W_T^{1,2}$ denote the subspace of $W_T^{1,2}$, consisting of functions of zero average.  To solve (\ref{8}) and (\ref{9}), we set up a map $(\eta,V) \ra (\xi,U) \, : R \times \overline W_T^{1,2} \ra R \times \overline W_T^{1,2}$, by solving 

\beqa
\lbl{10}
&  (\p (U'))'=-f(V)V'-k\sin \left( \eta+V(t) \right)+\frac{k}{T} \int_0^T \sin \left( \eta+V(t) \right) \, dt+e(t) \\ \nonumber
& \frac{k}{T} \int_0^T \sin \left( \xi+U(t) \right) \, dt=\mu \,. \nonumber
\eeqa
Since the right hand side of the first equation has average zero, this equation has  the unique $T$-periodic solution, by Lemma \ref{lma:2}. We claim that one can  find $\xi \in (-\frac{\pi}{2}, \frac{\pi}{2})$, solving the second equation in (\ref{10}).
Any solution of (\ref{10}) satisfies $|U'(t)|<a$ for all $x$, which implies that $\left( \int_0^T {U'(t)}^2 \, dt \right)^{1/2} \leq a \sqrt T$. Then by Lemma \ref{lma:3},  and our assumptions
\[
||U||_{L^{\infty}(R)} \leq \frac{\sqrt T}{2 \sqrt 3} \left( \int_0^T {U'(t)}^2 \, dt \right)^{1/2} \leq \frac{a  T}{2 \sqrt 3}<\frac{\pi}{2} \,.
\]
Then
\[
\sin(-\frac{\pi}{2}+U)<\sin(-\frac{\pi}{2}+\frac{a  T}{2 \sqrt 3})=-\cos \frac{a  T}{2 \sqrt 3}<0 \,,
\]
\[
\sin(\frac{\pi}{2}+U)>\sin(\frac{\pi}{2}+\frac{a  T}{2 \sqrt 3})=\cos \frac{a  T}{2 \sqrt 3}>0 \,.
\]
Hence, for any $\mu \in (-k\cos \frac{a  T}{2 \sqrt 3},k\cos \frac{a  T}{2 \sqrt 3})$, we can find a  $\xi \in (-\frac{\pi}{2}, \frac{\pi}{2})$, solving the second equation in (\ref{10}). Using Lemma \ref{lma:3.3}, we conclude that the map $(\eta,V) \ra (\xi,U)$ is a compact map of $R \times \overline W_T^{1,2}$ into $(-\frac{\pi}{2},\frac{\pi}{2}) \times B$, where $B$ is a  ball of radius $c_0$ in $\overline W_T^{1,2}$ ($c_0$ as in Lemma \ref{lma:3.3}). 
This map is also continuous. Indeed, by Lemma \ref{lma:2}, $U$ is continuous in $(\eta,V)$, and then writing the second equation in (\ref{10}) in the form $a \sin \xi+b \cos \xi=\mu$, with the coefficients $a$ and $b$ continuous in $(\eta,V)$, we see that $\xi$ continuous in $(\eta,V)$.
By Schauder's fixed point theorem there exists a fixed point, with $\xi \in  (-\frac{\pi}{2},\frac{\pi}{2})$. Similarly, $(\frac{\pi}{2},\frac{3\pi}{2}) \times B$ is   mapped into itself continuously and compactly, giving us a second fixed point, with $\xi \in (\frac{\pi}{2},\frac{3\pi}{2})$.
\epf

\section{Continuation of solutions}

We consider the following  linear periodic problem in the class of functions of zero average:  find a constant $\mu$, and $w(t) \in C^2(R)$ solving  
\beq
\lbl{11} \hspace{0.3in}
(a(t)w'(t))'+\la w'(t)+kh(t)w(t)=\mu, \s w(t+T)=w(t), \s \int_0^T w(s) \, ds=0,
\eeq
where $a(t) \in C^1(R)$ and $h(t) \in C(R)$ are  given  functions of period $T$, while $k$ and $\la$ are parameters. We denote by $C_T^2$ the subset of $C^2(R)$, consisting of $T$-periodic functions, and by $\bar C_T^2$ the subset of $C_T^2$, consisting of functions of zero average. Recall that $\omega =\frac{2 \pi}{T}$.

\begin{lma}\lbl{lma:5}
Assume that the $T$-periodic functions $a(t)$ and $h(t)$ satisfy $|h(t)| \leq 1$, $a(t)>a_0$ for all $t$,  and that  
\beq
\lbl{12}
k < a_0\omega \,.
\eeq
Then the only $T$-periodic solution of (\ref{11}) is $\mu=0$ and $w(t) \equiv 0$.
\end{lma}

\pf
Multiplying the equation in (\ref{11}) by $w(t)$ and integrating, we have (see Lemma \ref{lma:4})
\beqa \nonumber
& a_0 \int_0^T {w'}^2 \, dt \leq \int_0^T a(t) {w'}^2 \, dt=k\int_0^T h(t)ww'\, dt \\ \nonumber
& \leq k \left( \int_0^T w^2 \, dt \right)^{1/2} \left( \int_0^T {w'}^2 \, dt \right)^{1/2} \leq \frac{k}{\omega }  \int_0^T {w'}^2 \, dt \,, \nonumber
\eeqa
which contradicts (\ref{12}), unless $w(t) \equiv 0$, and then $\mu=0$.
\epf

\smallskip

We consider  the periodic problem 
\beqa
\lbl{14}
& (\p (u'))'+\la u'+ kg(u)=h(t)=\mu + e(t), \\
& u(t+T)=u(t), \s  u'(t+T)=u'(t)\,, \nonumber
\eeqa
which includes the case of relativistic pendulum with friction. The function $\p : (-a,a) \ra R$ is  assumed to be increasing, of class $C^2$,  $\lim _{u \ra \pm a} \p(u)=\pm \infty$, and  $\p (0)=0$.  The function $h(t) \in C(R)$ is assumed to be $T$-periodic, $g(u) \in C^1(R)$. We decompose $h(t)=\mu + e(t)$, with $\mu = \frac1T \int_0^T h(t) \, dt$, and $\int_0^T e(t) \, dt=0$. Decompose the solution $u(t)=\xi+U(t)$, with $\xi = \frac1T \int_0^T u(t) \, dt$, and $\int_0^T U(t) \, dt=0$. Integrating  (\ref{14}), 
\beq
\lbl{15}
\mu =\frac{k}{T} \int_0^T  g\left( \xi+U(t) \right) \, dt \,.
\eeq
Using this formula in (\ref{14}), we see that $U$ satisfies 

\beq
\lbl{16}
\s\s\s\s F(U) \equiv (\p (U'))'+\la U'+kg \left( \xi+U(t) \right)-\frac{k}{T} \int_0^T g \left( \xi+U(t) \right) \, dt=e(t) \,.
\eeq

Observe that $F: \bar C^2_T \ra \bar C_T$.

\begin{thm}\lbl{thm:5}
Assume that there is a constant $a_0$ such that
\beq
\lbl{17}
\p '(t) \geq a_0 \s \mbox{for $t \in (-a,a)$} \,,
\eeq
and that 
\beq
\lbl{18}
|g'(u) | \leq 1 \s \mbox{for all $u \in R$} \,.
\eeq
Assume finally that
\beq
\lbl{19}
k<a_0 \omega \,.
\eeq
Then  all solutions of (\ref{14}) lie on a unique continuous solution curve $(u,\mu)(\xi)$, with $\xi \in (-\infty,\infty)$.  Moreover, for any $\xi \in (-\infty,\infty)$ there exists a unique solution pair $(u(t),\mu)$ of (\ref{14}), with the average of $u(t)$ equal to $\xi$. I.e., $\xi$ is a global parameter on this solution curve.
\end{thm}

\pf
For each fixed $\xi$, finding the pair $(u,\mu)$ solving (\ref{14}) breaks down to first solving (\ref{16}) for $U$, and then finding $\mu $ from (\ref{15}). We show that the Implicit Function Theorem applies to (\ref{16}). The linearized operator is
\[
F'(U)V \equiv (\p' (U')V')'+\la V'+kg' \left( \xi+U(t) \right)V-\mu^* \,,
\]
where the constant $\mu^*$ stands for $\frac{k}{T} \int_0^T g' \left( \xi+U(t) \right) V(t) \, dt$. The operator $F'(U)V: \bar C^2_T \ra \bar C_T$ is injective because the only solution of 
\[
F'(U)V =0
\]
is $V=0$ and $\mu^* =0$, in view of the Lemma \ref{lma:5} and our conditions. By the Fredholm alternative this operator is also surjective, and hence the Implicit Function Theorem applies. This will allow us to continue solutions in the parameters $k$ and $\xi$.
\medskip

Turning to the existence and uniqueness of solutions, we embed  our problem into a family of problems
\beqa
\lbl{20}
& (\p (u'))'+\la u'+ \kappa g(u)=h(t)=\mu + e(t) \\ \nonumber
& u(t+T)=u(t), \s  u'(t+T)=u'(t)\\ \nonumber
& \frac1T \int_0^T u(t) \, dt=\xi \,, \nonumber
\eeqa
with $0 \leq \kappa \leq k$.  When $\kappa =0$ and $\mu =0$, the problem has a unique $T$-periodic solution of average $\xi$, by Lemma \ref{lma:2}. We now continue this solution in $\kappa$, i.e., we solve (\ref{20}) for $(u,\mu)$ as a function of $\kappa$ (keeping $\xi$ fixed). Again, we decompose the solution $u(t)=\xi+U(t)$, with $\xi = \frac1T \int_0^T u(t) \, dt$, and $\int_0^T U(t) \, dt=0$, and the Lyapunov-Schmidt decomposition (\ref{15}) and (\ref{16}) becomes
\beq
\lbl{21}
\mu =\frac{\kappa}{T} \int_0^T  g\left( \xi+U(t) \right) \, dt \,,
\eeq
\beq
\lbl{22}
\s\s\s\s F(U,\kappa) \equiv (\p (U'))'+\la U'+\kappa g \left( \xi+U(t) \right)-\frac{\kappa}{T} \int_0^T g \left( \xi+U(t) \right) \, dt=e(t) \,.
\eeq
The Implicit Function Theorem allows us to continue the solutions locally in $\kappa $ (first solving (\ref{22}) for $U$, and then finding $\mu $ from (\ref{21})). Since by Lemmas \ref{lma:3.3} and  \ref{lma:2.3}, solutions stay bounded, we can do the continuation for all $0 \leq \kappa \leq k$, obtaining the solution curve $(u,\mu)(\kappa)$ of  (\ref{14}), with the average of $u$ equal to $\xi$, and at $\kappa=k$, we have the desired solution. If we had another solution of average $\xi$, we would continue it for decreasing $\kappa$, obtaining a second solution of average $\xi$ at $\kappa=0$, in contradiction to Lemma \ref{lma:2}.
\medskip

Once we have a solution of (\ref{14}) at some $\xi$, we continue it in $\xi$, for all $-\infty< \xi<\infty$, by using the Implicit Function Theorem, as above.
\epf

This theorem implies that the curve $\mu =\mu (\xi)$ gives a faithful description of the existence and multiplicity of $T$-periodic solutions for the problem (\ref{14}). Properties of this curve can be described in more detail under further assumptions on $g(u)$. We begin with a result of Landesman-Lazer type.

\begin{thm}\lbl{thm:6}
Assume that the conditions of the Theorem \ref{thm:5} hold, and in addition, the function $g(u)$ has finite limits at $\pm \infty$, and
\[
g(-\infty)<g(u)<g(\infty) \s \mbox{for all $u \in R$} \,.
\]
Then the problem (\ref{14}) has a $T$-periodic solution if and only if 
\[
kg(-\infty)<\mu<kg(\infty) \,.
\]
\end{thm}

\pf
Necessity follows immediately from (\ref{15}). For sufficiency we also refer to (\ref{15}), and observe that we have a uniform bound on $U(t)$, when we do the continuation in $\xi$. Hence $\mu \ra kg(\pm \infty)$, as $\xi \ra \pm \infty$, and by continuity of $\mu (\xi)$, the problem (\ref{14}) is solvable for all $\mu$'s lying between these limits.
\epf

This result also follows from the Theorem $2$ in C. Bereanu and J. Mawhin \cite{BM}.
\medskip

\noindent
{\bf Example} We have solved the problem (\ref{14}) with $p(t)=\frac{t}{\sqrt{1-t^2}}$, $g(t)=\arctan t$, $e(t)=0.3 \sin \frac{2 \pi}{T} t$,   $T = 0.3$,  $k = 0.25$,  $\la = 0$.
The curve $\mu=\mu (\xi)$ is given in Figure $1$. We see that $\mu(\xi) \ra \pm \frac{\pi}{8}$, as $\xi \ra \pm \infty$. The picture also suggests the uniqueness of solutions.

\begin{figure}
\begin{center}
\scalebox{0.9}{\includegraphics{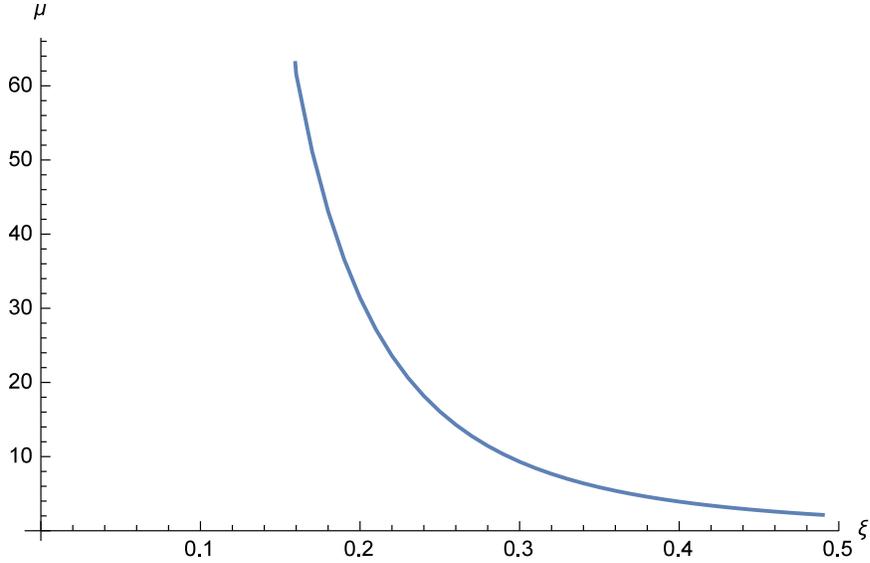}}
\end{center}
\caption{ An example for Theorem \ref{thm:6}}
\end{figure}
\medskip

A simple variation is provided by the following result, whose proof is similar.
\begin{thm}\lbl{thm:7}
Assume that the conditions of the Theorem \ref{thm:5} hold, and in addition,
\[
ug(u) >0 \s \mbox{for $|u|$ large}, \s \mbox{ and } \; \lim _{u \ra \pm \infty} g(u)=0 \,.
\]
Then there are constants $\mu _- <0< \mu _+$ so that the problem (\ref{14}) has at least two  $T$-periodic solution for $\mu \in (\mu _-,\mu _+) \setminus 0$, it has at least one  $T$-periodic solution for $\mu=\mu _- $, $\mu=0$ and $\mu=\mu _+ $, and no $T$-periodic solutions for $\mu$ lying outside of $(\mu _-,\mu _+)$.
\end{thm}

\pf
Since we have a uniform bound on $U(t)$, when we do the continuation in $\xi$, it follows from (\ref{15}) that for $\xi$ positive (negative) and large, $\mu$ is positive (negative) and it tends to zero as $\xi \ra \infty$ ($\xi \ra -\infty$).
\epf

\noindent
{\bf Example} We have solved the problem (\ref{14}) with $p(t)=\frac{t}{\sqrt{1-t^2}}$, $g(t)=\frac{3t}{1+t^2}$, $e(t)=0.45 \sin \frac{2 \pi}{T} t$,   $T = 0.2$,  $k = 0.1$,  $\la = 0.05$.
The curve $\mu=\mu (\xi)$ is given in Figure $2$. Here $\mu _- \approx -0.15$, and $\mu _+ \approx 0.15$.

\begin{figure}
\begin{center}
\scalebox{0.9}{\includegraphics{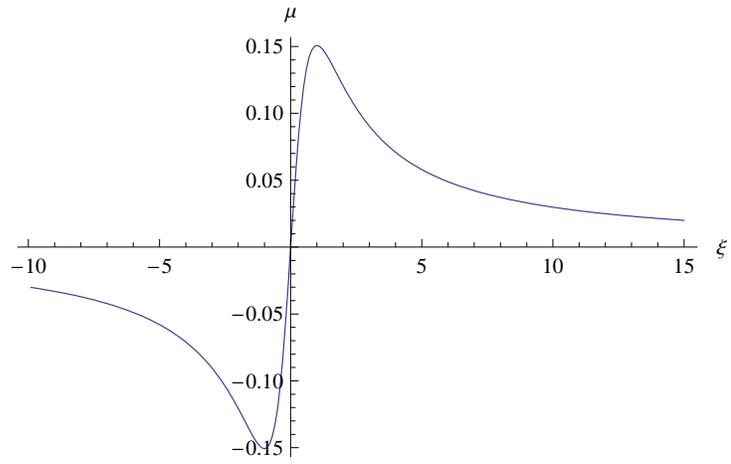}}
\end{center}
\caption{ An example for Theorem \ref{thm:7}}
\end{figure}

\begin{thm}\lbl{thm:8}
Assume that the conditions of the Theorem \ref{thm:5} hold, and in addition, the function $g(u)$ is periodic, of period $p$. Then the function $\mu =\mu(\xi)$ is periodic, of period $p$.
\end{thm}

\pf
The equation (\ref{16}) at $\xi +p$ is identical to the same equation at $\xi$, and hence $U(t)$ is the same. Then from (\ref{15}), $\mu (\xi +p)=\mu (\xi)$.
\epf

\noindent
{\bf Example} ({\em Relativistic pendulum}) We have solved the problem (\ref{14}) with $p(t)=\frac{t}{\sqrt{1-t^2}}$, $g(t)=\sin t$, $e(t)=0.15 \cos \frac{2 \pi}{T} t$,   $T =1$,  $k = 0.1$,  $\la = 0.1$.
The curve $\mu=\mu (\xi)$ is given in Figure $3$. The function $\mu =\mu(\xi)$ has period $2\pi$.

\begin{figure}
\begin{center}
\scalebox{0.9}{\includegraphics{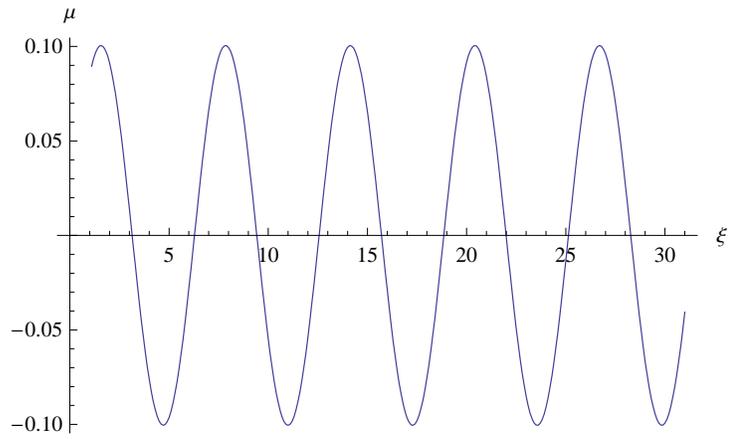}}
\end{center}
\caption{ An example for Theorem \ref{thm:8}}
\end{figure}
Each point in Figure $3$ represents a solution of the problem (\ref{14}), with period $1$. For example, when $\xi=5$, we have calculated $\mu \approx -0.09637$, and the actual $1$-periodic solution $u(t)$ is given in Figure $4$. We see very small variation of the solution $u(t)$ around its average value $\xi=5$. The variation is around $0.004$, compared with the variation of $0.15$ for the forcing term. We saw small variations for all values of the parameters, for all problems that we tried.

\begin{figure}
\begin{center}
\scalebox{0.9}{\includegraphics{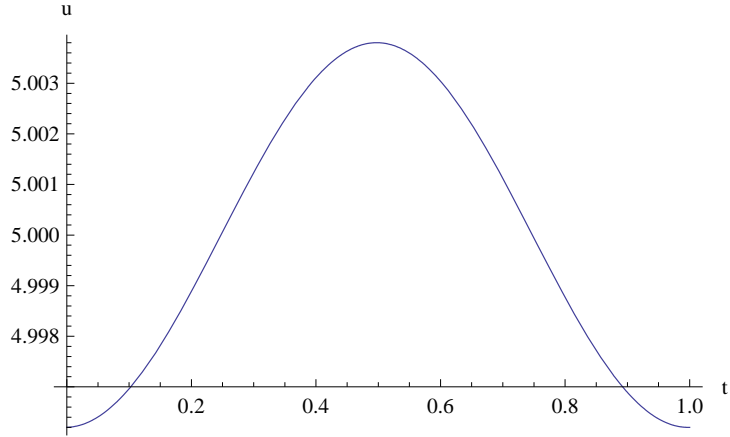}}
\end{center}
\caption{ A $1$-periodic solution of the problem (\ref{14})}
\end{figure}

We also have the following result of D.G. de Figueiredo  and   W.-M. Ni's \cite{FN} type, which does not restrict the behavior of $g(u)$ at infinity.

\begin{thm}
Assume that the conditions of the Theorem \ref{thm:5} hold, and in addition, the function $g(u)$ satisfies
\[
ug(u)>0 \s\s \mbox{for all $u \in R$} \,.
\]
Then for $\mu=0$, the problem (\ref{14}) has a solution.
\end{thm}

\pf
Proceeding as above, we have $\mu (\xi)>0$ ($<0$), when $|\xi|$ is large and $\xi$ is positive (negative). By continuity, $\mu (\xi_0)=0$ at some $\xi_0$.
\epf

This result also follows from the Theorem $2$ in C. Bereanu and J. Mawhin \cite{BM}.
\medskip

\section{Numerical computation of solutions}
\setcounter{equation}{0}
\setcounter{thm}{0}
\setcounter{lma}{0}

To find solutions of the periodic problem (\ref{14}), we used continuation in $\xi$, the average value of solution, and the Lyapunov-Schmidt decomposition given by (\ref{15}) and (\ref{16}), as described in the preceding section. We begin by implementing the numerical solution of the following periodic problem: given $T$-periodic functions $a(t)$, $b(t)$ and $f(t)$, and a constant $\la$, find the $T$-periodic solution of
\beq
\lbl{40}
\s\s\s\s L[y] \equiv \left(a(t) y'(t) \right)'+\la y' +b(t)y=f(t), \; y(t)=y(t+T), \; y'(t)=y'(t+T) \,.
\eeq
That turned out to be very simple, giving the capabilities of the {\em Mathematica} software. The general solution of (\ref{40}) is of course
\[
y(t)=Y(t)+c_1 y_1(t)+c_2 y_2(t) \,,
\]
where $Y(t)$ is a particular solution, and $y_1$, $y_2$ are two solutions of the corresponding homogeneous equation. To find $Y(t)$, we used the NDSolve command to solve (\ref{40}) with $y(0)=0$, $y'(0)=1$. {\em Mathematica} not only solves differential equations numerically, but it returns the solution as an interpolated function of $t$, practically indistinguishable from an explicitly defined function. We calculated $y_1$ and  $y_2$ by solving the corresponding homogeneous equation with the initial conditions $y_1(0)=0$, $y_1'(0)=1$ and $y_2(0)=1$, $y_2'(0)=0$. We then select $c_1$ and $c_2$, so that
\[
y(0)=y(T), \s y'(0)=y'(T) \,,
\]
which is just a linear $2 \times 2$ system. This gives us the $T$-periodic solution of (\ref{40}), or $L^{-1}[f(t)]$, where $L[y]$ denotes the left hand side of (\ref{40}), subject to the periodic boundary conditions. 
\medskip

Then we have implemented the ``linear solver", i.e., the numerical solution of the following problem: given $T$-periodic functions $a(t)$, $b(t)$ and $f(t)$, and a constant $\la$, find the constant $\mu^*$, so that the problem
\[
\left(a(t) y'(t) \right)'+\la y' +b(t)y=\mu^*+f(t), \; \int_0^T y(t) \, dt=0 
\]
has a $T$-periodic solution of zero average, and compute that solution $y(t)$. The solution is 
\[
y(t)=L^{-1}[f(t)]+\mu L^{-1}[1] \,,
\]
with the constant  $\mu$ chosen so that $\int_0^T y(t) \, dt=0$.
\medskip

Turning to the problem (\ref{14}), we begin with an initial $\xi _0$, and using a step size $\Delta \xi$, we compute solutions of average $\xi _i=\xi _0 +i \Delta \xi$, $i=1,2, \cdots, nsteps$, in the form $u=\xi _i+U$, where $U$ is the solution of  (\ref{16}) at $\xi =\xi _i$. Once $U$ is computed, we use (\ref{15}) to compute $\mu =\mu _i$. Finally, we plot the points $(\xi _i, \mu _i)$ to obtain the solution curve.
\medskip

To solve for $U(t)$, we apply Newton's method. If the iterate $U_n(t)$ is already computed, to solve for $U_{n+1}(t)$ we linearize the equation (\ref{16}) at $U_n(t)$, i.e., we apply the linear solver to find the $T$-periodic solution of
\[
\left(a(t)U_{n+1}'(t) \right)'+\la U_{n+1}' +b(t)U_{n+1}=\mu^*+f(t), \; \int_0^T U_{n+1}(t) \, dt=0 \,,
\]
with $a(t)= \p '(U_n') $, $b(t)=kg'(\xi _i+U_n)$, and $f(t)=\frac{d}{dt} \left( \p '(U_n') U_n'-\p (U_n') \right)+kg'(\xi _i+U_n)U_n-kg(\xi _i+U_n)+e(t)$.
The constant $\mu^*$ stands for 
\[
\mu^*=\frac{k}{T} \int_0^T \left[ g \left( \xi+U_n (t) \right)+g' \left( \xi+U_n (t) \right) \left( U_{n+1}(t) -U_n(t) \right) \right]\, dt \,.
\]
We found that two iterations of Newton's method, coupled with a relatively small $\Delta \xi$ (e.g., $\Delta \xi=0.1$), were sufficient for accurate computation of the solution curves.
\medskip

We have verified our numerical results by an independent calculation. Once  a periodic solution is computed at some $\mu$, we took its $u(0)$ and $u'(0)$, and computed numerically the solution with this initial data (using the NDSolve command). We had a complete agreement for all $\mu$, and all equations that we tried.
\medskip

\noindent
{\bf Acknowledgments} 
I wish to thank the referees for useful comments, and relevant references.

\end{document}